\documentclass[10pt,a4paper]{article}
\usepackage{amsmath,amssymb}
\usepackage[dvips]{graphicx}
\newcommand{\R}{{\mathbb R}} 
\newcommand{\Rn}{{\mathbb R}^n} 
\newcommand{\N}{{\mathbb N}} 

\newcommand{\irn}{\int\limits_{\Rn}}

\newcommand{\Wn}{W^{1,2}(\Rn)}

\newcommand{\dw}{\dot{w}}
\newcommand{\dz}{\dot{z}}
 
\newcommand{\wz}{\widetilde z}

\renewcommand{\a }{\alpha }
\renewcommand{\b }{\beta }
\renewcommand{\d }{\delta }
\newcommand{\D }{\Delta }
\newcommand{\e }{\varepsilon }

\newcommand{\n }{\nabla }

\newcommand{\s }{\sigma }

\newenvironment{pf}{\noindent{\sc Proof}.\enspace}{\rule{2mm}{2mm}\medskip}
\newtheorem{Theorem}{Theorem}[section]

\newtheorem{Lemma}[Theorem]{Lemma}

\newtheorem{remark}[Theorem]{Remark}
\newtheorem{remarks}[Theorem]{Remarks}
\newenvironment{Remark}{\begin{remark} \rm}{\rule{2mm}{2mm}\end{remark}}
\author{ A. Ambrosetti \and A. Malchiodi \and S. Secchi }

\title{Multiplicity results for some nonlinear \\
Schr\"{o}dinger equations with potentials}

\date{Preprint SISSA  108/2000/AF}
\begin{document}
\maketitle

\section{Introduction}\label{sec:int}

This  paper deals with the existence of multiple positive solutions for a class
of
nonlinear Schr\"{o}dinger Equation (NLS in short)
\begin{equation}\tag{NLS}\label{eq:00}
\left\{
\begin{array}{lll} -\e^2 \D u + u +V(x)u=K(x)u^p, && \\
u \in W^{1,2}(\Rn),\,u>0 &&
\end{array}
\right.
\end{equation}
where $\D$ denotes the Laplace operator and
$$
1<p<2^*=\left\{ \begin{array}{lll} &\frac{2n}{n-2},&\text{ if } \,
n\geq 3, \\
 &+\infty, & \text{ if } \, n=2.
\end{array}
\right.
$$
In the sequel we will always assume that $V,K:\Rn\to \R$ satisfy
\begin{description}
\item[$(V1)$] $V\in C^{2}(\Rn)$, $V$ and $D^{2}V$ are bounded;

 \item[$(V2)$] $\inf[1+V(x)]>0$;

\item[$(K1)$] $K\in C^{2}(\Rn)$, $K$ is bounded
and $K(x)>0$ $\forall\;x\in \Rn$.
\end{description}
One seeks solutions $u_\e$ of (NLS) that concentrate, as $\e \to 0$, near some
point $x_0\in \Rn$ (semiclassical standing waves).  By this we mean that
for all $x\in \Rn\setminus\{x_{0}\}$ one has that $u_\e(x)\to 0$ as $\e\to 0$.

 When $K$ equals a positive constant, say $K(x)\equiv 1$,
 (NLS) has been widely investigated, see
 \cite{ABC,ABer,dPF,FW,Li,Oh,W} and references therein.
 Moreover, the existence of multibump solutions has also been
studied in \cite{CN,G}, see also
\cite{ABer} where solutions with infinitely many bumps has been proved.
 It has been also pointed out,
 see e.g. \cite[Section 6]{ABC}, that the results contained in the
forementioned papers
 can be  extended to equations where $u^p$ is substituted by a
function $g(u)$, which behaves like $u^p$. Nonlinearities depending upon
$x$ have been
handled in \cite{Gr,WZ} where the existence of one-bump solutions is proved.

In a group of papers Equation (NLS) is studied
by perturbation arguments. For example, in \cite{ABC} a Liapunov-Schmidt
type procedure
is used to
reduce, for $\e$ small, (NLS) to a finite dimensional equation, that
inherits the
variational structure of the original problem. So, one looks for the
critical points of a finite
dimensional functional, which leading term is strictly related to the
behaviour of $V$ near its stationary points.

A different approach has been carried out in another group of papers, like
\cite{dPF,W,WZ}. First, one finds a solution $u$ of (NLS) as a Mountain
Pass critical point of
the Euler functional;  after, one proves the concentration as $\e \to 0$.
Although this second approach allows to deal with potentials $V$ such that
$\sup V =+\infty$, it
works, roughly speaking, near minima of $V$, only. Actually, as pointed out
in \cite{ABC}, see also \cite[Theorem 3.2]{AB}, the Morse index of $u$ is
equal to
$1\,+$ the Morse index of the stationary point $V$  where concentration
takes place.
Hence $u$ can be found as a Mountain Pass only in the case of minima of $V$.
Such a severe restriction is even more apparent when one deals with
nonlinearities depending
upon $x$, see \cite{WZ}.  A case dealing with saddle-like potentials
is studied in \cite{dPFM}, but under some technical conditions
involving the saddle level and the supremum of $V$.

The second approach sketched above has been used in \cite{CL1} to obtain
a multiplicity result.
Letting $V_0=\inf [1+V(x)]$, it is assumed that
$K(x)\equiv 1$ and $\liminf_{|x|\to\infty}(1+V(x))>V_0 >0$. It has been
shown that, for $\e>0$
small,
(NLS) has at least $cat(M,M_\d)$ solutions, where
 $M=\{x\in \Rn:1+V(x)=V_0\}$, $M_\d$ is a $\d$ neighbourhood of $M$ and
 $cat$ denotes the Lusternik-Schnirelman category. This result has been extended
 in \cite{CL2} to  nonlinearities like $K(x)u^p + Q(x)u^q$, with
 $1<q<p$.

 In the the present paper we will improve the above multiplicity result.
 For example, when $K(x)\equiv 1$, we can obtain the results sumarized below.

 \begin{description}
 \item[$(i)$]  suppose that $V$ has a manifold $M$ of stationary points,
and  that $M$
 is nondegenerate in the sense of Bott, see \cite{Bott}. We prove, see
Theorem \ref{th:main},
 that (NLS) has at least $l(M)$ critical points concentrating near points
of $M$. Here
 $l(M)$ denotes the cup long of $M$.  This kind of result is new because
 \cite{CL1,CL2} deal only with (absolute) minima.

 \item[$(ii)$] If the points of $M$ are local minima or maxima of $V$
 the above result can be sharpened because $M$ does not need to be a
 manifold and $l(M)$ can be subsituted by $cat(M,M_\d)$, see Theorem
 \ref{th:CL}.  Unlike \cite{CL1,CL2}, we can handle sets of {\it local}
 minima.  Furthermore, we do not require any condition at infinity.
 The case of maxima is new.

 \item[$(iii)$] When $K(x)$ is not identically constant, the preceding result
 holds true provided that $V$ is substituted by
 $A=(1+V)^{\theta}K^{-2/(p-1)}$, where $\theta=-(n/2)+(p+1)/(p-1)$, see
 Theorem \ref{th:K}.
 \end{description}
A similar multiplicity result holds for problems involving more
general nonlinearities, like in \cite{Gr}, see Remark \ref{rem:CLW}.

 Our proof relies on perturbation arguments, variational in nature,
as in \cite{AB,ABC}.  Neverhtless, a straight application of the
arguments used in those papers, would only provide the existence of one
solution of (NLS) because that procedure leads to find critical points
of a finite dimensional functional which does not inherit the
topological features of $M$.  So, in order to find multiplicity
results like those described above, it is necessary to use a different
finite dimensional reduction.  With this new approach one looks for
critical points of a finite dimensional functional $\Phi_{\e}$ which
is defined in a tubular neighbourhood $\cal B$ of $M$ and, roughly, is
close to $\n V$ (or $\n A$) on its boundary $\partial {\cal B}$.  It
turns out that $\cal B$ is an isolating block for the flow of $\n \Phi$.
Then one can use the Conley theory \cite{Con}, in particular some
results of \cite{CZ}, to find find $l(M)$ solutions.  This new
procedure is carried out in Sections \ref{eq:inv} and \ref{sec:fdr},
while the main existence results are proved in Section \ref{sec:main}.

\

{\bf Notation}
\begin{itemize}
\item $\Wn$ denotes the usual Sobolev space, endowed with the norm
$$
\|u\|^{2}=\int_{\Rn}|\nabla u|^{2}dx + \irn u^2 dx.
$$

	\item  $o_h (1)$ denotes a function that tends to $0$ as $h\to 0$.

	\item  $c, C, c_i$ etc. denote (possibly different) positive constants,
	independent of the parameters.
\end{itemize}

\section{Preliminaries}\label{sec:prel}
In this section we will collect some preliminary material which will be
used in the rest of the paper.
 In order to semplify the notation, we will discuss in detail
equation (NLS) when $K\equiv 1$.  The general case will be handled in
Theorem \ref{th:K} and requires some minor changes, only.

 Without loss of generality we can assume that $V(0)=0$.
Performing the change of variable $x\mapsto \e x$, equation $(NLS)$ becomes
\begin{equation}\tag{$P_{\e}$}
 - \D u + u +V(\e x)u=u^p.
\label{eq:P}
\end{equation}
Solutions of (\ref{eq:P}) are the critical points $u\in \Wn$ of
 $$
 f_\e (u)=f_0(u)+ \frac{1}{2}\int_{\Rn}V(\e x)u^2dx,
 $$
 where
$$
f_0(u)=\frac{1}{2}\|u\|^2 - \frac{1}{p+1}\int_{\Rn}u^{p+1}dx.
$$
The solutions of (\ref{eq:P}) will be found near a solution of
\begin{equation}
- \D u + u +V(\e \xi)u=u^p,
\label{eq:xi}
\end{equation}
for an appropriate choice of $\xi\in \Rn$.  The solutions of
(\ref{eq:xi}) are critical points
of the functional
\begin{equation}
F^{\e\xi}(u)=f_{0}(u)+\frac{1}{2}\,V(\e
\xi)\,\int_{\Rn}u^2dx-\frac{1}{p+1}\int_{\Rn}u^{p+1}dx
\label{eq:F}
\end{equation}
and can be found explicitly.  Let $U$ denote the unique, positive,
radial solution of
\begin{equation}\label{eq:unp}
 - \D u + u =u^p ,\qquad u \in W^{1,2}(\Rn).
\end{equation}
Then a straight calculation shows that $\a U(\b x)$ solves (\ref{eq:P})
whenever
$$
\b=\b(\e\xi)= [1+V(\e\xi)]^{1/2}\quad {\rm and}\quad \a = \a(\e\xi)=[\b
(\e\xi)]^{2/(p-1)}.
$$
We set
\begin{equation}
	z^{\e\xi}(x)=\a(\e\xi)U(\b(\e\xi)x)
\label{eq:zU}
\end{equation}
and
$$
Z^{\e}=\{z^{\e\xi}(x-\xi):\xi\in \Rn\}.
$$
 When there is no possible misunderstanding, we will write $z$, resp. $Z$,
 instead of $z^{\e\xi}$, resp $Z^{\e}$.
 We will also use the notation $z_{\xi}$ to denote the function
 $z_{\xi}(x):=z^{\e\xi}(x-\xi)$.
   Obviously all the
 functions in $z_{\xi}\in Z$ are solutions of (\ref{eq:unp}) or, equivalently,
 critical points of $F^{\e\xi}$.
 For future references let us point out some estimates. First of all, we
evaluate:
\begin{eqnarray*}
&& \partial_{\xi}z^{\e\xi}(x-\xi)= \partial_{\xi}
\left[\a(\e\xi)U(\b(\e\xi)(x-\xi))\right] \\
&& =  \e\a'(\e\xi)U(\b(\e\xi)(x-\xi))+
 \e\a(\e\xi)\b'(\e\xi)U'(\b(\e\xi)(x-\xi))-\a(\e\xi)U'(\b(\e\xi)(x-\xi)).
\end{eqnarray*}
Recalling the
 definition of $\a$, $\b$ one finds:
 \begin{equation}
 \partial_{\xi}z^{\e\xi}(x-\xi)=-\partial_{x}z^{\e\xi}(x-\xi)+O(\e|\nabla
V(\e\xi)|).	
 	\label{eq:partialz}
 \end{equation}

\noindent The next Lemma shows that $\nabla f_{\e}(z_\xi)$ is close to
zero when
$\e$ is small.

 \begin{Lemma}\label{lem:1}
For all $\xi\in \Rn$ and all $\e>0$ small, one has that
$$
\|\nabla f_{\e}(z_{\xi})\|\leq C\left(\e |\nabla
V(\e\xi)|+\e^{2}\right),\quad C>0.
$$
 \end{Lemma}
\begin{pf}
From
$$
	f_{\e}(u)  =  F^{\e\xi}(u)+\frac{1}{2}
	\int_{\Rn}\left[V(\e x)-V(\e\xi)\right]u^{2}dx
$$
and since $z_{\xi}$ is a critical point of $F^{\e\xi}$, one has
\begin{eqnarray*}
(\nabla f_{\e}(z_{\xi})|v)	 & = & (\nabla F^{\e\xi}(z_{\xi})|v)+
\int_{\Rn}\left[V(\e
x)-V(\e\xi)\right]z_{\xi} v\,dx  \\
& = & \int_{\Rn}\left[V(\e x)-V(\e\xi)\right]z_{\xi} v\,dx . \\
\end{eqnarray*}
Using the H\"{o}lder inequality, one finds
\begin{equation}
|(\nabla f_{\e}(z_{\xi})|v)|^{2}\leq \|v\|^{2}_{L^{2}} \int_{\Rn}|V(\e
x)-V(\e\xi)|^{2}z_{\xi}^{2}dx.
\label{eq:1.1}
\end{equation}
From the assumption that $|D^{2}V(x)|\leq {\rm const.}$ one infers
$$
|V(\e x)-V(\e\xi)|\leq \e |\nabla V(\e\xi)|\cdot |x-\xi|+c_{1}\e^{2}|x-\xi|^{2}.
$$
This implies
\begin{eqnarray}
 &  &\int_{\Rn}|V(\e x)-V(\e\xi)|^{2}z_{\xi}^{2}dx\leq  \nonumber \\
 & & c_1 \leq\e^{2}|\nabla V(\e\xi)|^{2}\int_{\Rn}|x-\xi|^{2}z^{2}(x-\xi)dx +
 c_{2}\e^{4}\int_{\Rn}|x-\xi|^{4}z^{2}(x-\xi)dx.
\label{eq:1.3}
\end{eqnarray}
Recalling (\ref{eq:zU}), a direct calculation yields
\begin{eqnarray*}
\int_{\Rn}|x-\xi|^{2}z^{2}(x-\xi)dx & = &
\a^{2}(\e\xi)\int_{\Rn}|y|^{2}U^{2}(\b(\e\xi) y)dy \\
& = & \a^{2} \b^{-n-2}\int_{\Rn}|y'|^{2}U^{2}(y')dy'\leq c_{3}.
\end{eqnarray*}
From this (and a similar calculation for for the last integral in the
above formula) one infers
\begin{equation}
\int_{\Rn}|V(\e x)-V(\e\xi)|^{2}z_{\xi}^{2}dx\leq c_{4}\e^{2}|\nabla
V(\e\xi)|^{2} + c_{5}\e^{4}.
\label{eq:1.4}
\end{equation}
Putting together (\ref{eq:1.4}) and (\ref{eq:1.1}), the Lemma follows.
 \end{pf}

\section{Invertibility of $D^{2}f_{\e}$ on $TZ^{\perp}$}\label{sec:inv}
In this section we will show that $D^{2}f_{\e}$ is invertible on
$TZ^{\perp}$.  This will be the main tool to perform the finite
dimensional reduction, carried out in Section \ref{sec:fdr}.

Let $L_{\e,\xi}:(T_{z_{\xi}}Z^{\e})^{\perp}\to
(T_{z_{\xi}}Z^{\e})^{\perp}$ denote the operator defined by setting
$(L_{\e,\xi}v|w)= D^{2}f_{\e}(z_{\xi})[v,w]$.  We want to show
\begin{Lemma}\label{lem:inv}
Given $\overline{\xi}>0$ there exists $C>0$ such that for $\e$ small enough
one has that
\begin{equation}
	|(L_{\e,\xi}v|v)|\geq C \|v\|^{2},\qquad \forall\;|\xi|\leq
\overline{\xi},\;\forall\; v\in
	(T_{z_{\xi}}Z^{\e})^{\perp}.
\label{eq:inv}
\end{equation}
\end{Lemma}
\begin{pf}
From (\ref{eq:partialz}) it follows that every element
$\zeta\in T_{z_{\xi}}Z$ can be written in the form $\zeta =
z_{\xi}-\partial_{x}z^{\e\xi}(x-\xi)+O(\e)$.  As a consequence,
\begin{description}
	\item[$(*)$]  it suffices to prove (\ref{eq:inv}) for all
$v\in {\rm span}\{z_{\xi},\phi\}$, where $\phi$ is orthogonal to
 ${\rm span}\{z_{\xi},\partial_{x}z^{\e\xi}(x-\xi)\}$.
\end{description}
Precisely we shall prove that there exist $C_{1},C_{2}>0$ such that
for all $\e>0$ small and all $|\xi|\leq |\overline{\xi}|$
one  has:
\begin{eqnarray}
	(L_{\e,\xi}z_{\xi}|z_{\xi})& \leq & - C_{1}< 0.
	\label{eq:neg} \\
(L_{\e,\xi}\phi|\phi)&\geq & C_{2} \|\phi\|^2.
	\label{eq:claim}
\end{eqnarray}
It is clear that the Lemma immediately follows from $(*)$,
(\ref{eq:neg}) and (\ref{eq:claim}).

\

\noindent {\bf Proof of} (\ref{eq:neg}). First, let us recall that, since
$z_{\xi}$ is a Mountain Pass critical point of $F$, then given
$\overline{\xi}$ there exists $c_0>0$ such that for all $\e>0$ small
and all $|\xi|\leq |\overline{\xi}|$ one finds:
\begin{equation}
	 D^2 F^{\e\xi}(z_{\xi})[z_{\xi},z_{\xi}] < -c_0< 0.
	\label{eq:D2Fa}
\end{equation}
One has:
$$
	(L_{\e,\xi}z_{\xi}|z_{\xi})=D^2 F^{\e\xi}(z_{\xi})[z_{\xi},z_{\xi}]+
	\int_{\Rn}\left[V(\e x)-V(\e\xi)\right]z_{\xi}^2 dx.
$$
The last integral can be estimated as in (\ref{eq:1.4}) yielding
\begin{equation}
	(L_{\e,\xi}z_{\xi}|z_{\xi})\leq D^{2} F^{\e\xi}(z_{\xi})[z_{\xi},z_{\xi}]
	+ c_1 \e |\nabla V(\e\xi)|+c_2 \e^2.
	\label{eq:D2zz}
\end{equation}
From (\ref{eq:D2Fa}) and (\ref{eq:D2zz}) it follows that
(\ref{eq:neg}) holds.

\

\noindent {\bf Proof of} (\ref{eq:claim}).  As before, the fact that
$z_{\xi}$ is a Mountain Pass critical point of $F$ implies that
\begin{equation}
 D^2 F^{\e\xi}(z_{\xi})[\phi,\phi]>c_1 \|\phi\|^2.
	\label{eq:D2Fb}
\end{equation}
Let $R\gg 1$ and consider a radial smooth function
$\chi_{1}:\R^{n}\mapsto \R$ such that
\begin{equation}\tag{$\chi'$}\label{eq:c1}
\chi_{1}(x) = 1, \quad \hbox{ for } |x| \leq R; \qquad
\chi_{1}(x) = 0, \quad \hbox{ for } |x| \geq 2 R;
\end{equation}
\begin{equation}\tag{$\chi''$}\label{eq:c2}
|\n \chi_{1}(x)| \leq \frac{2}{R}, \quad \hbox{ for } R \leq |x| \leq 2 R.
\end{equation}
We also set $ \chi_{2}(x)=1-\chi_{1}(x)$.
Given $\phi$ let us consider the functions
$$
\phi_{i}(x)=\chi_{i}(x-\xi)\phi(x),\quad i=1,2.
$$
A straight computation yields:
$$
\irn \phi^2 = \irn \phi_1^2 + \irn \phi_2^2 + 2\irn \phi_{1} \,
\phi_{2},
$$
$$
\irn |\n \phi|^2 = \irn |\n \phi_1|^2 + \irn |\n \phi_2|^2 + 2\irn
\n\phi_{1} \cdot \n \phi_{2},
$$
and hence
$$
\| \phi \|^2 = \| \phi_1 \|^2 + \| \phi_2 \|^2+
2 \irn\left[ \phi_{1} \, \phi_{2}+\n\phi_{1} \cdot \n
\phi_{2}\right].
$$
Letting $I$ denote the last integral,
one immediately finds:
$$
I=\underbrace{\irn \chi_{1}\chi_{2}(\phi^{2}+|\n \phi|^{2})}_{I_{\phi}} +
\underbrace{\irn\phi^{2}\n\chi_{1}\cdot \n\chi_{2}}_{I'}+
\underbrace{\irn\phi_{1}\n\phi\cdot\n\chi_{2}+\phi_{2}\n
\phi\cdot\n\chi_{1}}_{I''}.
$$
Due to the definition of $\chi$, the two integrals $I'$ and $I''$
reduce to integrals from $R$ and $2R$, and thus they are $o_{R}(1)\|\phi\|^{2}$.
As a consequence we have that
\begin{equation}\label{eq:d}
\| \phi \|^2 = \| \phi_1 \|^2 + \| \phi_2 \|^2 + 2I_\phi + o_R(1)\| \phi \|^2,
\end{equation}
After these preliminaries, let us evaluate the three terms in the
equation below:
$$
(L_{\e,\xi}\phi|\phi)=
\underbrace{(L_{\e,\xi}\phi_{1}|\phi_{1})}_{\a_{1}}+
\underbrace{(L_{\e,\xi}\phi_{2}|\phi_{2})}_{\a_{2}}+
2\underbrace{(L_{\e,\xi}\phi_{1}|\phi_{2})}_{\a_{3}}.
$$
One has:
\begin{equation} \label{eq:alfa1}
\a_{1}=(L_{\e,\xi}\phi_{1}|\phi_{1})=D^{2}F^{\e\xi}[\phi_{1},\phi_{1}]+\irn\left[V(\e x)-V(\e\xi)\right]\phi_{1}^{2}.
\end{equation}
In order to use (\ref{eq:D2Fb}), we introduce the function
$\overline{\phi}_{1}=\phi_{1}-\psi$, where
$$
\psi=(\phi_{1}|z_{\xi})z_{\xi}+(\phi_{1}|\partial_{x}z_{\xi})\partial_{x}z_{
\xi}.
$$
Then we have:
\begin{equation}
D^{2}F^{\e\xi}[\phi_{1},\phi_{1}]=D^{2}F^{\e\xi}[\overline{\phi}_{1},\overline{\phi}_{1}]+
D^{2}F^{\e\xi}[\psi,\psi]+2D^{2}F^{\e\xi}[\overline{\phi}_{1},\psi]
\label{eq:alfa2}
\end{equation}
 Let us explicitely point out that $\overline{\phi}_{1}\perp {\rm
span}\{z_{\xi},\partial_{x}z^{\e\xi}(x-\xi)\}$ and hence
(\ref{eq:D2Fb}) implies
\begin{equation}
D^{2}F^{\e\xi}[\overline{\phi}_{1},\overline{\phi}_{1}]\geq c_{1}
\|\overline{\phi}_{1}\|^{2}.
\label{eq:alfa3}
\end{equation}
On the other side, since $(\phi|z_{\xi})=0$ it follows:
\begin{eqnarray*}
	(\phi_{1}|z_{\xi}) & = & (\phi|z_{\xi})-(\phi_{2}|z_{\xi})=
-(\phi_{2}|z_{\xi}) \\
& = & -\irn\phi_{2}z_{\xi}dx-\irn \n z_{\xi}\cdot \n \phi_{2}dx \\
&=& -\irn\chi_{2}(y)z(y)\phi(y+\xi)dy-\irn \n z(y)\cdot \n
\chi_{2}(y)\phi(y+\xi)dy .
\end{eqnarray*}
Since $\chi_{2}(x)=0$ for all $|x|<R$, and since $z(x)\to 0$ as
$|x|=R\to\infty$, we infer
$(\phi_{1}|z_{\xi})=o_{R}(1)\|\phi\|$.  Similarly one shows that
$(\phi_{1}|\partial_{x}z_{\xi})=o_{R}(1)\|\phi\|$ and it follows that
\begin{equation}
	 \|\psi\|=o_{R}(1)\|\phi\|.
\label{eq:alfa4}
\end{equation}
We are now in position to estimate the last two terms in Eq. (\ref{eq:alfa2}).
Actually, using (\ref{eq:alfa4}) we get
\begin{equation}
D^{2}F^{\e\xi}[\psi,\psi]=\|\psi\|^{2}+V(\e\xi)\irn\psi^{2}-p\irn
z_{\xi}^{p-1}\psi^{2}=o_{R}(1)\|\phi\|^{2}.	
\label{eq:alfa5}
\end{equation}
The same arguments readily imply
\begin{equation}
D^{2}F^{\e\xi}[\overline{\phi}_{1},\psi]=(\overline{\phi}_{1}|\psi)
+V(\e\xi)\irn \overline{\phi}_{1}\psi -p\irn
z_{\xi}^{p-1}\overline{\phi}_{1}\psi =o_{R}(1)\|\phi\|^{2}.
\label{eq:alfa6}
\end{equation}
Putting together (\ref{eq:alfa3}), (\ref{eq:alfa5}) and  (\ref{eq:alfa6})
we infer
\begin{equation}
D^{2}F^{\e\xi}[\phi_{1},\phi_{1}]\geq \|\phi_{1}\|^{2}+o_{R}(1)\|\phi\|^{2}.
\label{eq:alfa7}
\end{equation}
Using arguments already carried out before, one has
\begin{eqnarray*}
\irn|V(\e x)-V(\e\xi)|\phi_{1}^{2}dx & \leq
&c_{2}\irn|x-\xi|\chi^{2}(x-\xi)\phi^{2}(x)dx \\
& \leq & \e c_{3}\irn y\chi(y)\phi^{2}(y+\xi)dy \\
& \leq & \e c_{4}\|\phi\|^{2}.
\end{eqnarray*}
This and (\ref{eq:alfa7}) yield
\begin{equation}
	\a_{1}=(L_{\e,\xi}\phi_{1}|\phi_{1})\geq c_{5}\|\phi_{1}\|^{2}-
	\e c_{4}\|\phi\|^{2}+o_{R}(1)\|\phi\|^{2}.
\label{eq:alfa8}
\end{equation}
Let us now estimate $\a_{2}$.  One finds
$$
\a_{2}=(L_{\e,\xi}\phi_{2}|\phi_{2})=\irn |\n \phi_{2}|^{2}+
\irn (1+V(\e x))\phi_{2}^{2} -p\irn z_{\xi}^{p-1}\phi_{2}^{2}
$$
and therefore, using $(V2)$,
$$
\a_{2}  \geq  c_{6} \|\phi_{2}\|^{2} -p\irn z_{\xi}^{p-1}\phi_{2}^{2}.
 $$
As before, $\phi_{2}(x)=0$ for all $|x|<R$ and $z(x)\to 0$ as
$|x|=R\to\infty$ imply that
\begin{equation}
\a_{2}  \geq  c_{6} \|\phi_{2}\|^{2}+o_{R}(1)\|\phi\|^{2}.
\label{eq:alfa9}
\end{equation}
In a quite similar way one shows that
\begin{equation}
\a_{3} \geq c_{7}I_{\phi}+o_{R}(1)\|\phi\|^{2}.
\label{eq:alfa10}
\end{equation}
Finally, (\ref{eq:alfa8}), (\ref{eq:alfa9}), (\ref{eq:alfa10}) and the fact
that $I_{\phi}\geq 0$, yield
\begin{eqnarray*}
(L_{\e,\xi}\phi|\phi) & = & \a_{1}+\a_{2}+2\a_{3}\\
& \geq & c_{8}\left[\|\phi_{1}\|^{2}+\|\phi_{2}\|^{2}
+2 I_{\phi}\right]-c_{9}\e\|\phi\|^{2'}+ o_{R}(1)\|\phi\|^{2}.
\end{eqnarray*}
Recalling (\ref{eq:d}) we infer that
$$
(L_{\e,\xi}\phi|\phi)\geq c_{10}\|\phi\|^{2}-c_{9}\e\|\phi\|^{2'}+
o_{R}(1)\|\phi\|^{2}.
$$
Taking $\e$ small and $R$ large, eq.  (\ref{eq:claim}) follows.
This completes the proof of Lemma \ref{lem:inv}.
\end{pf}

\section{The finite dimensional reduction}\label{sec:fdr}
 In this Section we will show that the existence of critical points of $f_{\e}$
 can be reduced to the search of critical points of an auxiliary finite
 dimensional functional.  The proof will be carried out in 2
 subsections dealing, respectively, with  a Liapunov-Schmidt reduction, and
 with the behaviour of the auxiliary finite dimensional functional.
In a final subsection we handle the general case in which $K$ is not
identically constant.

\subsection{A Liapunov-Schmidt type reduction}\label{subsec:LS} The main
result of this section is the following lemma.

\begin{Lemma}\label{lem:w}
For $\e>0$ small and $|\xi|\leq \overline{\xi}$ there exists a unique
$w=w(\e,\xi)\in
(T_{z_\xi} Z)^{\perp}$ such that
$\nabla f_\e (z_\xi + w)\in T_{z_\xi} Z$.
Such a $w(\e,\xi)$
 is of class $C^{2}$, resp.  $C^{1,p-1}$, with respect to $\xi$, provided that
 $p\geq 2$, resp.  $1<p<2$.
Moreover, the functional $\Phi_\e (\xi)=f_\e (z_\xi+w(\e,\xi))$ has
the same regularity of $w$ and satisfies:
 $$
\n \Phi_\e(\xi_0)=0\quad \Longrightarrow\quad \n
f_\e\left(z_{\xi_0}+w(\e,\xi_0)\right)=0.
$$
\end{Lemma}
\begin{pf}
Let $P=P_{\e\xi}$ denote the projection onto $(T_{z_\xi} Z)^\perp$. We want
to find a
solution $w\in (T_{z_\xi} Z)^{\perp}$ of the equation
$P\nabla f_\e(z_\xi +w)=0$.  One has that $\n f_\e(z+w)=
\n f_\e (z)+D^2 f_\e(z)[w]+R(z,w)$ with $\|R(z,w)\|=o(\|w\|)$, uniformly
with respect to
$z=z_{\xi}$, for  $|\xi|\leq \overline{\xi}$. Using the notation introduced
in the
preceding Section \ref{sec:inv}, we are led to the equation:
$$
L_{\e,\xi}w + P\n f_\e (z)+PR(z,w)=0.
$$
According to Lemma \ref{lem:inv}, this is equivalent to
$$
w = N_{\e,\xi}(w), \quad \mbox{where}\quad
N_{\e,\xi}(w)=-L_{\e,\xi}^{-1}\left( P\n f_\e (z)+PR(z,w)\right).
$$
From Lemma \ref{lem:1} it follows that

\begin{equation}	\label{eq:N}
	\|N_{\e,\xi}(w)\|\leq c_1 (\e|\n V(\e\xi)|+\e^2)+ o(\|w\|).
	\end{equation}
Then one readily checks that $N_{\e,\xi}$ is a contraction on some ball in
$(T_{z_\xi} Z)^{\perp}$
provided that $\e>0$ is small enough and $|\xi|\leq \overline{\xi}$.
Then there exists a unique $w$ such that $w=N_{\e,\xi}(w)$.  Let us
point out that we cannot use the Implicit Function Theorem to find
$w(\e,\xi)$, because the map $(\e,u)\mapsto P\n f_\e (u)$ fails to be
$C^2$.  However, fixed $\e>0$ small, we can apply the Implicit
Function Theorem to the map $(\xi,w)\mapsto P\n f_\e (z_\xi + w)$.
Then, in particular, the function $w(\e,\xi)$ turns out to be of class
$C^1$ with respect to $\xi$.  Finally, it is a standard argument, see
\cite{AB,ABC}, to check that the critical points of $\Phi_\e
(\xi)=f_\e (z+w)$ give rise to critical points of $f_\e$.
\end{pf}

\begin{Remark}\label{rem:w}
 From (\ref{eq:N}) it immediately follows that:
\begin{equation}	\label{eq:w}
	\|w\|\leq C \left(\e |\n V(\e\xi)|+\e^2\right),
\end{equation}
where $C>0$.
 \end{Remark}

\

\noindent  For future references, it is
convenient to estimate the derivative $\partial_\xi w$.

\begin{Lemma}\label{lem:Dw}
One has that:
\begin{equation}	\label{eq:Dw}
	\|\partial_\xi w\|\leq c \left(\e |\n
V(\e\xi)|+O(\e^2)\right)^\gamma,\qquad c>0,\;
	\gamma=\min\{1,p-1\}.
\end{equation}
\end{Lemma}
\begin{pf}
In the proof we will write $\dw_i$, resp $\dz_i$, to denote the components of
 $\partial_\xi w$, resp.
$\partial_\xi z$; moreover we will set $h(z,w)=(z+w)^p-z^p-pz^{p-1}w$.
With these notations, and recalling that
$L_{\e,\xi}w = -\D w +w +V(\e x)w -pz^{p-1}w$, it follows that
$w$ satisfies $\forall \;v\in (T_{z_\xi} Z)^{\perp}$:
$$
(w|v) + \irn V(\e x)wvdx- p\irn z^{p-1}wvdx +\irn [V(\e x)-V(\e\xi)]zv dx-
\irn h(z,w)v dx=0.
$$
Hence $\dw_i$ verifies
\begin{eqnarray}
	 &  & (\dw_i |v) + \irn V(\e x)\dw_i vdx- p\irn z^{p-1}\dw_i vdx
-p(p-1)\irn z^{p-2}\dz_i wv dx
	 \nonumber  \\
	 &  & +\irn [V(\e x)-V(\e\xi)]\dz_i v dx-\e \partial_{x_i} V(\e\xi)\irn zv dx
	  - \irn (h_z\dz_i +
h_w \dw_i )v dx=0.\label{eq:dw}
\end{eqnarray}
Let us set $L'=L_{\e,\xi}-h_w $. Then
 (\ref{eq:dw}) can be written as
\begin{equation}\label{eq:Ltilde}
(L'\dw_i|v) =p(p-1)\irn z^{p-2}\dz_i wv
- \irn [V(\e x)-V(\e\xi)]\dz_i v
+\e \partial_{x_i} V(\e\xi) \irn zv + \irn  h_z\dz_i v ,
	\end{equation}
and hence one has
$$
|(L'\dw_i|v)| \leq c_1 \|w\|\cdot\|v\|+\left|\irn [V(\e x)-V(\e\xi)]\dz v
dx\right|
+c_2 \e |\n V(\e\xi)|\cdot\|v\|
+ \left|\irn h_z\dz v dx\right|.
$$
As in the proof of Lemma \ref{lem:1} one gets
$$
\left|\irn [V(\e x)-V(\e\xi)]\dz v dx\right|\leq\left( c_3 \e |\n
V(\e\xi)|+c_4 \e^2\right)\,\|v\|.
$$
Furthermore, the definition of $h$ immediately yields
$$
\left|\irn h_z\dz v dx\right|\leq c_5 \|w\|^\gamma \|v\|,\quad \mbox{where}\quad
\gamma = \min\{1,p-1\}.
$$
Putting together the above estimates we find
$$
|(L'\dw_i|v)|\leq
\left[c_6 \e |\n V(\e\xi)|+c_4 \e^2 + c_6 \|w\|^\gamma\right]\,\|v\|
$$
Since $h_w\to 0$ as $w \to 0$, the operator
$L'$, likewise $L$, is  invertible for $\e >0$ small and therefore one finds
$$
\|\dw\|\leq \left( c_7 \e |\n V(\e\xi)|+c_8 \e^2\right) + c_9 \|w\|^\gamma,
$$
Finally, using (\ref{eq:w})
the Lemma follows.	
 \end{pf}

\subsection{The finite dimensional functional}\label{subs:3}
The main purpose of this subsection is to use the estimates
on $w$ and $\partial_\xi w$ estabilished above to find
an expansion of $\n \Phi_\e (\xi)$, where $\Phi_\e (\xi)= f_\e(z_\xi
+w(\e,\xi))$.  In
the sequel, to be short, we will often write $z$ instead of $z_\xi$
and $w$ instead of $w(\e,\xi)$.  It is always understood that
$\e$ is taken in such a way that all the results discussed in the
preceding Section hold.

For the reader's convenience we will divide the arguments in some steps.

\

\noindent {\sl Step 1.} We have:
$$
	\Phi_\e (\xi) =  \frac{1}{2}\|z+w\|^2 + \frac{1}{2}\irn V(\e x)(z+w)^2
	-\frac{1}{p+1}\irn (z+w)^{p+1}.
$$
Since $-\D z+z+V(\e\xi)z=z^p$ we infer that
\begin{eqnarray*}
	\|z\|^2 & = &  -V(\e\xi)\irn z^2 + \irn z^{p+1}, \\
	(z|w) & = & -V(\e\xi)\irn zw + \irn z^pw,
\end{eqnarray*}
Then we find :
\begin{eqnarray*}
	 \Phi_\e (\xi)& = & \left(\frac{1}{2}-\frac{1}{p+1}\right) \irn z^{p+1} +
	 \frac{1}{2}\int_{\Rn}\left[V(\e x)-V(\e\xi)\right]z^{2}\\
	 &  & +\int_{\Rn}\left[V(\e x)-V(\e\xi)\right]zw+
	 \frac{1}{2}\irn V(\e x)w^2 \\
&&+ \frac{1}{2} \|w\|^2 - \frac{1}{p+1}\irn
\left[(z+w)^{p+1}-z^{p+1}-(p+1)z^pw\right].
\end{eqnarray*}
Since $z(x)=\a(\e\xi)U(\b(\e\xi)x)$, where
$\a=(1+V)^{1/(p-1)}$ and $\b=(1+V)^{1/2}$, see (\ref{eq:zU}),
it follows
$$
\irn z^{p+1}dx= C_0 \left(1+V(\e\xi)\right)^\theta,\quad C_0=
 \irn U^{p+1};\qquad \theta=\frac{p+1}{p-1}-\frac{n}{2}.
$$
Letting $C_1= C_0 [1/2 -1/(p+1)]$ one has
\begin{eqnarray}
	\Phi_\e (\xi)  =  C_1 \left(1+V(\e\xi)\right)^\theta+
	\frac{1}{2}\int_{\Rn}\left[V(\e x)-V(\e\xi)\right]z^{2}+
	\int_{\Rn}\left[V(\e x)-V(\e\xi)\right]zw\nonumber && \\
	 + \frac{1}{2}\irn V(\e x)w^2 +\frac{1}{2} \|w\|^2 -
	\frac{1}{p+1} \irn \left[(z+w)^{p+1}-z^{p+1}-(p+1)z^pw\right].&&
	 \label{eq:F1}
\end{eqnarray}

\

\noindent {\sl Step 2.}
Let us now evaluate the derivative of the right hand side
of (\ref{eq:F1}).
For this, let us write:

\begin{equation}	\label{eq:D}
\Phi_\e (\xi)=C_1(1+V(\e\xi))^{\theta}+\Lambda_\e(\xi)+
\Psi_\e(\xi),
	\end{equation}
where
$$
\Lambda_\e(\xi)=\frac{1}{2}\int_{\Rn}\left[V(\e x)-V(\e\xi)\right]z^2+
\int_{\Rn}\left[V(\e x)-V(\e\xi)\right]zw
$$
and
$$
\Psi_\e(\xi)=
	\frac{1}{2}\irn V(\e x)w^2 +\frac{1}{2} \|w\|^2 -
	 \irn \left[(z+w)^{p+1}-z^{p+1}-(p+1)z^pw\right].
$$
It is easy to check, by a direct calculation, that
\begin{equation}	\label{eq:D1}
	|\n \Psi_\e(\xi)|\leq c_1\|w\|\,\|\dw\|.
	\end{equation}
Furthermore,
since $V(\e x)-V(\e\xi)=\e\n V(\e\xi)\cdot(x-\xi)+\e^2 D^2V(\e\xi+\tau
\e(x-\xi))[x-\xi,x-\xi]$ for some $\tau\in [0,1]$ we get:
\begin{eqnarray*}
	\int_{\Rn}\left[V(\e x)-V(\e\xi)\right]z^2dx & = & \e \irn \n
V(\e\xi)\cdot(x-\xi)z^2dx\\
&&\quad +\e^2 \irn D^2V(\e\xi+\tau
\e(x-\xi))[x-\xi,x-\xi]z^2 dx  \\
&=&  \e \irn \n V(\e\xi)\cdot y\, z^2(y)dy\\
	 &  &\quad +\e^2 \irn D^2V(\e\xi+\tau\e y)[y,y]z^2(y) dy \\
&=& \e^2 \irn D^2V(\e\xi+\tau\e y)[y,y]z^2(y) dy.
\end{eqnarray*}
Similarly, from $V(\e x)-V(\e\xi)=\e\n V(\e\xi+\tau
\e(x-\xi))\cdot (x-\xi)$ one finds
$$
\int_{\Rn}\left[V(\e x)-V(\e\xi)\right]zw=\e \irn \n V(\e\xi+\tau
\e(x-\xi))\cdot (x-\xi) zw
$$
The preceding equations imply:
\begin{equation}
	|\n \Lambda_\e(\xi)|
	\leq c_2 \e^2 + c_3 \e\|w\|.
	\label{eq:D2}
\end{equation}
Taking the gradients in (\ref{eq:D}), using (\ref{eq:D1}-\ref{eq:D2}) and
recalling
the estimates (\ref{eq:w}) and (\ref{eq:Dw}) on $w$ and $\dw$, respectively,
we readily find:
\begin{Lemma}\label{lem:DF}
Let $a(\e\xi)=\theta C_1(1+V(\e\xi))^{\theta -1}$. Then one has:
$$
\n \Phi_\e (\xi)= a(\e\xi) \e\n V(\e\xi) + \e^{1+\gamma} R_\e(\xi),
$$
where $|R_\e(\xi)|\leq const$ and $\gamma=\min\{1,p-1\}$.
\end{Lemma}

\begin{Remark}\label{rem:Phi}
Using similar arguments one can show that
 $$
 \Phi_\e (\xi)= C_1 \left(1+V(\e\xi)\right)^\theta + \rho_\e(\xi),\quad
C_1>0,\quad
 \theta=\frac{p+1}{p-1}-\frac{n}{2},
 $$
 where $|\rho_\e(\xi)|\leq const. \left(\e |\nabla V(\e\xi)| +\e^2\right)$.
\end{Remark}

\subsection{The general case}\label{subsec:gen}
Let us indicate the counterpart of the above results in the case in which
$K$ is not constant.  Since the arguments are quite similar, we
will only  outline the main modifications that are needed.

\

After rescaling, the solutions of (NLS) are the critical points of
$$
\widetilde{f}_\e(u)=\frac{1}{2}\|u\|^2 +\frac{1}{2}\int_{\Rn}V(\e x)u^2-
\frac{1}{p+1}\int_{\Rn}K(\e x)u^{p+1}.
$$
The solutions of (NLS) will be found near  solutions of
\begin{equation}
- \D u + u +V(\e \xi)u=K(\e \xi)u^p,
\label{eq:xiK}
\end{equation}
namely near critical points of
$$
\widetilde{F}^{\e\xi}(u)=\frac{1}{2}\|u\|^2
+\frac{1}{2}\,V(\e \xi)\,\int_{\Rn}u^2 -
\frac{1}{p+1}\,K(\e \xi)\int_{\Rn}u^{p+1}.
$$
If $\wz=\wz^{\e\xi}$ is a solution of (\ref{eq:xiK}), then $\wz (x)=
\tilde{\a}(\e\xi)U(\tilde{\b}(\e\xi)x)$, where
$$
\tilde{\a}(\e\xi)=\left(1+V(\e\xi)\right)^{1/2},\quad \tilde{\b}(\e\xi)=
\left(\frac{1+V(\e\xi)}{K(\e\xi)}\right)^{1/(p-1)}.
$$
This implies that

\begin{equation}	\label{eq:wzA}
\irn \wz ^{p+1}=C_0 A(\e\xi),
	\end{equation}
where
$$
A(x)=\left(1+V(x)\right)^{\theta}[K(x)]^{-2/(p-1)},\quad
\theta=\frac{p+1}{p-1}-\frac{n}{2}.
$$
Define $\widetilde{L}=\widetilde{L}_{\e,\xi}$ on $(T\widetilde{Z})^{\perp}$
by setting $(\widetilde{L} v|w)=D^2\widetilde{f}_\e(\wz_\xi)[v,w]$. As in Lemma
\ref{lem:inv} one shows that $\widetilde{L}$ is invertible for $\e$ small.
Furthermore, one has
$$
\widetilde{f}_\e(u)=\widetilde{F}^{\e\xi}(u)+
\frac{1}{2}\irn [V(\e x)-V(\e \xi)]u^2 -
\frac{1}{p+1}\irn[K(\e x)-K(\e \xi)]u^{p+1}
$$
and, as in Lemma \ref{lem:1} one finds that
$$
\|\widetilde{f}_\e(u)\|\leq c_1 \left(\e(|\n V(\e \xi)|+|\n
K(\e\xi)|)+\e^2\right).
$$
This and the invertibility of $\widetilde{L}$ imply,
as in  Lemma \ref{lem:w}, the existence of
$\widetilde{w}(\e,\xi)$ such that the critical points of the finite
dimensional functional
$$
\widetilde{\Phi}_\e(\xi)=\widetilde{f}_\e(\wz_\xi+\widetilde{w}(\e,\xi))
$$
give rise to critical points of $\widetilde{f}_\e(u)$.
Such a $\widetilde{w}$ is $C^1$ with respect
to $\xi$ and, as in Lemma \ref{lem:Dw}, the following estimate holds
$$
\|\partial_\xi\widetilde{w}\|\leq c_2
\left(\e(|\n V(\e \xi)|+|\n K(\e\xi)|)+O(\e^2)\right)^\gamma.
$$
It remains to study the finite dimensional functional
$\widetilde{\Phi}_\e (\xi)$,
$$
\widetilde{\Phi}_\e(\xi)=\frac{1}{2}\|\wz+\widetilde{w}\|^2 +
\frac{1}{2}\irn V(\e x)(\wz+\widetilde{w})^2
	-\frac{1}{p+1}\irn K(\e x)(\wz+\widetilde{w})^{p+1}.
$$
Since now $\|\wz\|^2=-V(\e\xi)\int \wz^2 + K(\e\xi)\int \wz^{p+1}$
and $(\wz |\widetilde{w})=-V(\e\xi)\int \wz\widetilde{w}+ K(\e\xi)\int \wz^{p}
\widetilde{w}$, one gets
\begin{eqnarray*}
	 \widetilde{\Phi}_\e (\xi)& = & \irn\left(\frac{1}{2}K(\e x)-\frac{1}{p+1}
	 K(\e\xi)\right) \wz^{p+1} +
	 \frac{1}{2}\irn\left[V(\e x)-V(\e\xi)\right]\wz^{2}\\
	 &  & +\int_{\Rn}\left[V(\e x)-V(\e\xi)\right]\wz \widetilde{w}+
	 \frac{1}{2}\irn V(\e x)\widetilde{w}^2+ \frac{1}{2} \|w\|^2 \\
	 && +K(\e\xi)\irn\wz^p \widetilde{w} -
	 \frac{1}{p+1}\irn K(\e x) \left[(z+w)^{p+1}-z^{p+1}-(p+1)z^pw\right]
\end{eqnarray*}
From $K(\e x)=K(\e\xi)+\e\n K(\e\xi)\cdot (x-\xi)+O(\e^2)$ and since
$\int \n K(\e\xi)\cdot y \wz^{p+1}=0$ we infer:
$$
\irn\left(\frac{1}{2}K(\e x)-\frac{1}{p+1}
	 K(\e\xi)\right) \wz^{p+1}=
	 (\frac{1}{2}-\frac{1}{p+1})\irn
	 K(\e\xi) \wz^{p+1} + O(\e^2).
 $$
Using (\ref{eq:wzA}), one finds the counterpart of (\ref{eq:F1}), namely
\begin{multline*}
	\widetilde{\Phi}_\e (\xi)  =  C_1 A(\e\xi)+
	\frac{1}{2}\irn\left[V(\e x)-V(\e\xi)\right]\wz^{2}
	  +\int_{\Rn}\left[V(\e x)-V(\e\xi)\right]\wz \widetilde{w}
	+ \frac{1}{2}\irn V(\e x)\widetilde{w}^2\\
	+ \frac{1}{2} \|w\|^2
 +K(\e\xi)\irn\wz^p \widetilde{w} - \frac{1}{p+1}\irn K(\e x)
 \left[(z+w)^{p+1}-z^{p+1}-(p+1)z^pw\right]+O(\e^2).
\end{multline*}
Taking the derivative of the above equation and using
the preceding estimates, one finally yields

\begin{equation}\label{eq:wDF}
\n \widetilde{\Phi}_\e (\xi)=C_1 \n A(\e\xi)+\e^{1+\gamma} \widetilde{R}_\e
(\xi),
\quad \gamma=\min\{1,p-1\},
 \end{equation} where $|\widetilde{R}_\e (\xi)|\leq const$, which is
the counterpart of Lemma \ref{lem:DF}.  Let us also point out that,
like in Remark \ref{rem:Phi}, one has:
\begin{equation}
	\widetilde{\Phi}_\e (\xi)=C_1 A(\e\xi)+O(\e).
	\label{eq:A}
\end{equation}

\section{Main results}\label{sec:main}
In this section we will prove the main results of the present paper.
First, some preliminaries are in order.

Given a set $M\subset \Rn$, $M\ne \emptyset$, we denote by $M_\d$ its $\d$
neighbourhood.
If $M\subset N$, $cat(M,N)$ denotes the Lusternik-Schnirelman category
of $M$ with respect to $N$, namely the least integer $k$ such that
$M$ can be covered by $k$ closed subsets of
$N$, contractible to a point in $N$. We set $cat(M)=cat(M,M)$.

The {\sl cup long} $l(M)$ of $M$ is defined by
$$
l(M)=1+\sup\{k\in \N: \exists\, \a_{1},\ldots,\a_{k}\in
\check{H}^{*}(M)\setminus 1, \,\a_{1}\cup\ldots\cup\a_{k}\ne 0\}.
$$
If no such class exists, we set $l(M)=1$. Here $\check{H}^{*}(M)$ is the
Alexander
cohomology of $M$ with real coefficients and $\cup$ denotes the cup product.

If $\Phi\in C^{2}(\Rn,\R)$, the
(steepest descent) gradient flow of $\Phi$ is the flow $\s^{t}(x)$ generated
by $ -\n \Phi$, namely the solutions of the Cauchy problem $d\s/dt=
-\n \Phi(\s),\;\s^{0}(x)=x$.
When $\Phi$ is merely $C^{1}$, $ -\n \Phi$ is replaced by
the {\sl psudogradient field} related to $\Phi$, see \cite{Pal}, and the
corresponding
flow will be called the {\sl pseudogradient flow} of $\Phi$.


%

 The following result is nothing but Theorems 4, 5 and the
subsequent Corollary in \cite{CZ}, adapted to our notations and
purposes.

\begin{Theorem} [\cite{CZ}] \label{th:CZ}
 Let $M\subset \Rn$ be a smooth compact manifold, $D_{1}=\{x_{1}\in
 \R^{d_{1}}: |x_{1}|\leq k_{1}\}$, $D_{2}=\{x_{2}\in \R^{d_{2}}:
 |x_{2}|\leq k_{2}\}$ and suppose that ${\cal B}=M\times D_{1}\times
 D_{2}$ is an isolating block (see \cite{Con} for the definition) for
 the gradient flow of a smooth functional $\Phi:\Rn\to \R$.  Then the
 gradient flow has an invariant set $S$ in $\cal B$ and $l(S)\geq
 l({\cal B})=l(M)$.  See \cite[Theorem 4]{CZ}.  Furthermore, $S$
 contains at least $l(S)$ critical points of $\Phi$.
\end{Theorem}
Roughly, $\cal B$ is a compact tubular neighbourhood of $M$,
$M\times D_1$ is the {\sl exit set}, namely the set of points in
$\partial\cal B$ where the $\n \Phi$ is inward pointing and $M\times D_2$ is
the entrance set, namely the set of points in $\partial\cal B$ where
the $\n \Phi$ is outward pointing.  A specific property of an
isolating block is that the exit set is closed relative to the set
$\{x\in {\cal B}:\exists \,t>0 \,\mbox{s.t.}\,\sigma^{t}(x)\not\in{\cal
B}\}$.  In particular, this implies that if $\cal B$ is an isolating
block for a flow, then it is also an isolating block for any flow
close the given one.  For more details we refer to \cite{Con}.

\begin{Remark}\label{rem:C}
It is easy to see that when $\Phi$ is of class $C^{1}$, Theorem
\ref{th:CZ} still holds, using the {\sl pseudo-gradient flow} of
$\Phi$ instead of the gradient flow.
\end{Remark}

We are now ready to state our multiplicity results.
Our first result deals with (NLS) with $K(x)\equiv 1$,
namely with the equation (\ref{eq:P}) introduced in Section
\ref{sec:prel}.  Let us suppose that $V$ has a smooth manifold of
critical points $M$.  According to Bott \cite{Bott}, we say that $M$
is nondegenerate (for $V$) if every $x\in M$ is a nondegenerate
critical point of $V_{|M^{\perp}}$.  The Morse index of $M$ is, by definition,
the Morse index of any $x\in M$, as
critical point of $V_{|M^{\perp}}$.
\begin{Theorem}\label{th:main}
Let  $(V1-2)$ hold and suppose $V$ has
a nondegenerate smooth manifold of critical points  $M$.
Then for $\e>0$ small, (\ref{eq:P}) has at least $l(M)$ solutions that
concentrate near
points of $M$.
\end{Theorem}
\begin{pf}First of all, we fix $\overline{\xi}$ in such a way that
$|x|<\overline{\xi}$ for all $x\in M$.  We will apply the
finite dimensional procedure with such $\overline{\xi}$ fixed. Since $M$ is
a nondegenerate smooth manifold of critical points of $V$, we can find
an isolating block ${\cal B}=M\times D_{1}\times D_{2}$ for the
gradient flow of $V$, where the {\sl exit set} $D^1$ has dimension
$k_{1}=Index(M)$.  As pointed out before Remark \ref{rem:C}, $\cal B$
is stable under small perturbation of the flow.  Then
from Lemma \ref{lem:DF}  we  deduce that
$$
{\cal B}^{\e}=\{\xi:\e\xi\in {\cal B}\}
$$
is also
an isolating block for $\Phi_{\e}$, provided $\e$ is small enough.  An
application of Theorem \ref{th:CZ} yields the existence of at least
$l(M)$ critical points of $\Phi_{\e}$.
Actually, when $1<p<2$ and
$\Phi_\e$ is merely of class $C^{1,p-1}$, we use Remark \ref{rem:C}.
If $\xi_{i}$ is any of those critical points, Lemma \ref{lem:w}
implies that $u_{\e,\xi_{i}}=z^{\e\xi_{i}}(x-\xi_{i})+w(\e,\xi_{i})$
is a critical point of $f_{\e}$.  Finally $u_{\e,\xi_{i}}(x/\e)$ is a
solution of (\ref{eq:P}) that concentrates near $\xi_{i}$.  This
completes the proof.
\end{pf}

The next result deals with the more general equation (NLS).
The results will be given using the auxiliary function
$$
A(x)=\left(1+V(x)\right)^{\theta}[K(x)]^{-2/(p-1)},\quad
\theta=\frac{p+1}{p-1}-\frac{n}{2},
$$
introduced in Subsection \ref{subsec:gen}.

\begin{Theorem}\label{th:K}
Let  $(V1-2)$ and $(K1)$ hold and suppose $A$ has
a nondegenerate smooth manifold of critical points  $\widetilde{M}$.
Then  for $\e>0$ small, (NLS) has at least $l(\widetilde{M})$ solutions
that concentrate near
points of $\widetilde{M}$.
\end{Theorem}
\begin{pf}
The proof is quite similar to the preceding one, by using the results
discussed in Subsection \ref{subsec:gen}.
\end{pf}
\begin{Remark}\label{rem:CLW}
The preceding results can be extended to cover a class of nonlinearities
$g(x,u)$ satisfying the same assumptions of \cite{Gr}.  In addition to
some techenical conditions, one roughly assumes that the problem

\begin{equation}
	-\D u+u+V(\e\xi)u=g(\e\xi,u),\quad u>0,\quad u\in \Wn,
\label{eq:Gr}
\end{equation}
has a unique radial solution $z=z^{\e\xi}$ such that the linearized problem
at $z$ is invertible on $(T_{z}Z)^{\perp}$.  In such a case one can
obtain the same results as above with the auxiliary function $A$
substituted by
$$
{\cal A}(\e\xi)=\irn\left[\frac{1}{2}g(\e\xi,z^{\e\xi}(x-\xi))-
G(\e\xi,z^{\e\xi}(x-\xi))\right]dx,
$$
where $\partial_{u}G=g$.  When $g(x,u)=K(x)u^p$ one finds that
${\cal A}=A$.  But, unlike such a case, the function $\cal A$ cannot be written
in an explicit way.  This is the reason why
 we have focused our study to the model problem (NLS)
when the auxiliary function $A$ has an explicit and neat form.  Let us
also point out that the class of nonlinearities handled in \cite{CL2}
does not even require that (\ref{eq:Gr}) has a unique solution.
\end{Remark}

\

\noindent When we deal with local minima (resp. maxima) of $V$, or $A$, the
preceding results
can be improved because the number of positive solutions of (NLS)
can be estimated by means of the category and $M$ does not need to be
a manifold.

\begin{Theorem}\label{th:CL}
Let  $(V1-2)$ and $(K1)$ hold and suppose $A$ has
a compact set $X$ where $A$ achieves a strict local minimum, resp.
maximum.

Then  there exists $\e_{\d}>0$ such that (NLS) has at least $cat(X,X_\d)$
solutions that concentrate near points of $X_{\d}$, provided $\e\in
(0,\e_{\d})$.
\end{Theorem}
\begin{pf}
Let $\d>0$ be such that
$$
b:=\inf\{A(x):x\in \partial X_{\d}\}>a:= A_{|X}
$$
and fix again
$\overline{\xi}$ in such a way that $X_{\d}$ is contained in
$\{x\in\Rn : |x|<\overline{\xi}\}$.  We set
$
X^{\e}=\{\xi:\e\xi\in X\}$, $X_{\d}^{\e}=\{\xi:\e\xi\in X_{\d}\}$ and
$Y^{\e}=\{\xi\in X_{\d}^{\e} :\widetilde{\Phi}_{\e}(\xi)\leq C_{1}(a+b)/2\}$.
From
(\ref{eq:A}) it follows that there exists $\e_{\d}>0$ such that
\begin{equation}
	X^{\e}\subset Y^{\e}\subset X^{\e}_{\d},
\label{eq:X}
\end{equation}
 provided $\e\in
(0,\e_{\d})$.  Moreover, if $\xi\in \partial X_{\d}^{\e}$ then
$V(\e\xi)\geq b$ and hence
$$
\widetilde{\Phi}_{\e}(\xi)\geq C_{1}V(\e\xi) + o_{\e}(1) \geq C_{1}b +
o_{\e}(1) .
$$
On the other side, if $\xi\in Y^{\e}$
then $\widetilde{\Phi}_{\e}(\xi)\leq C_{1}(a+b)/2$.  Hence, for
 $\e$ small, $Y^{\e}$ cannot meet $\partial X_{\d}^{\e}$
 and this readily implies that $Y^{\e}$ is compact.
Then $\widetilde{\Phi}_{\e}$ possesses at least $cat(Y^{\e},X^{\e}_{\d})$
critical points in $ X_{\d}$.  Using (\ref{eq:X}) and the properties of
the category one gets
$$
cat(Y^{\e},Y^{\e})\geq cat(X^{\e},X^{\e}_{\d})=cat(X,X_{\d}),
$$
and the result follows.
\end{pf}

\

\begin{Remark}
The approach carried out in the present paper
works also in the more standard case when $V$, or $A$, has an isolated
set $S$ of stationary points and $deg(\n V,\Omega,0)\not= 0$ for some
neighbourhood $\Omega$ of $S$.  In this way we can, for example, recover the
results of \cite{Gr}.
\end{Remark}

\

\begin{center}
{\bf Acknowledgment}
\end{center}
 The work has been supported by M.U.R.S.T.  under the national
project "Variational methods and nonlinear differential equations"

\end{document}